\documentclass[12pt,a4paper]{article}

\usepackage{preamble}

\title{Maps of Mori Dream Spaces in Cox coordinates\\
       Part I: existence of descriptions}
\author{\JaBu \and Oskar K{\k e}dzierski}

\newdir{smallhead}{
:(15,1)\dir^{>}
*:(15,-1)\dir_{>}}

\newdir{multihead}{
\dir{smallhead}
*!(0,1.1)\dir{smallhead}
*!(0,-1.1)\dir{smallhead}
*{\POS(-6,0); \POS(0, 1.1) **\crv{~**\dir{-} (-4.5,1.1)}}
*{\POS(-6,0); \POS(0,-1.1) **\crv{~**\dir{-} (-4.5,-1.1)}}
}

\newcommand{\multito}{
\
\begin{xy}
  =(0,0) "SOURCE";
  =(9,0) "TARGET";
  {\ar@{multihead} "SOURCE"; "TARGET"}
\end{xy}
\
}

\newcommand{\multimapsto}{
\
\begin{xy}
  =(0,0) "SOURCE";
  =(9,0) "TARGET";
  {\ar@{|-multihead} "SOURCE"; "TARGET"}
\end{xy}
\
}

\DeclareMathOperator{\Reg}{Reg}
\newcommand{\agree}[2]{{\operatorname{Agr}(#1, #2)}}
\newcommand{\irrel}[1]{\operatorname{Irrel}(#1)}
\newcommand{\Dind}{D_{\mathrm{ind}}}
\newcommand{\Dirrel}{D_{\mathrm{irrel}}}

\newcommand{\starfan}[1]{\operatorname{Star}(#1)}
\newcommand{\quofan}{\Sigma_{Z}}

\usepackage{color}
\newenvironment{red}{\color{red}}{}
\newcommand{\bred}{\begin{red}}
\newcommand{\ered}{\end{red}}

\newenvironment{blue}{\color{blue}}{}
\newcommand{\bblue}{\begin{blue}}
\newcommand{\eblue}{\end{blue}}

\newenvironment{green}{\color{green}}{}
\newcommand{\bgreen}{\begin{green}}
\newcommand{\egreen}{\end{green}}

\usepackage[textsize=tiny]{todonotes}
\begin{document}

\maketitle

\begin{abstract}
   Any rational map between affine spaces, projective spaces or toric varieties
      can be described in terms of their affine, homogeneous, or Cox coordinates.
   We show an analogous statement in the setting of Mori Dream Spaces.
   More precisely (in the case of regular maps) we show there exists a finite extension of the Cox ring of the source,
      such that the regular map lifts to a morphism from the Cox ring of the target to the finite extension.
   Moreover the extension only involves roots of homogeneous elements.
   Such a description of the map can be applied in practical computations.
\end{abstract}

\medskip
{\footnotesize
\noindent\textbf{addresses:} \\
J.~Buczy\'nski, \verb|jabu@mimuw.edu.pl|,
   Institute of Mathematics of the Polish Academy of Sciences, ul.~\'Sniadeckich 8, 00-656 Warsaw, Poland,
 and Faculty of Mathematics, Computer Science and Mechanics, University of Warsaw, ul.~Banacha 2, 02-097 Warszawa, Poland\\
O.~K{\k e}dzierski, \verb|oskar@mimuw.edu.pl|,
   Faculty of Mathematics, Computer Science and Mechanics, University of Warsaw, ul.~Banacha 2, 02-097 Warszawa, Poland

\noindent\textbf{keywords:}
Mori Dream Space, Cox ring, rational maps, coordinate description, multi-valued map.


\noindent\textbf{financial support:}
The authors are supported by Polish National Science Center (NCN), project 2013/11/D/ST1/02580. Buczy\'nski is also supported
          by a scholarship of Polish Ministry of Science.
This article was partially written during the Polish Algebraic Geometry mini-Semester (miniPAGES),
   which was supported by the grant 346300 for IMPAN from the Simons Foundation and the matching 2015-2019 Polish MNiSW fund.}

\section{Introduction}

By definition, an algebraic morphism of two affine varieties $\varphi\colon X\to Y $
  is a geometric interpretation of an algebra homomorphism $\varphi^* \colon B \to A$ of their affine coordinate rings.
Here $X = \Spec A$ and $Y = \Spec B$.
If $X = \PP^m$ and $Y = \PP^n$ instead,
  and $A \simeq \kk[\fromto{x_0}{x_m}]$ and $B \simeq \kk[\fromto{y_0}{y_n}]$ are their homogeneous coordinate rings,
  then any algebraic morphism $\varphi\colon \PP^m \to \PP^n$ is determined a homomorphism $B \to A$
  satisfying the usual homogeneity and base point freeness conditions.
Rational maps between affine varieties or projective spaces have similar interpretations in terms of the fields of
  fractions of coordinate rings.

Furthermore, an analogous coordinatewise description applies also to a rational map $\varphi\colon X \dashrightarrow Y$
   between toric varieties $X$ and $Y$ \cite{brown_jabu_maps_of_toric_varieties}, \cite{cox_functor}.
In this case the underlying algebra involves the \emph{Cox rings} of toric varieties, also called their \emph{homogeneous coordinate rings},
   or \emph{total coordinate rings}.
The major difference between previous settings is that we need to extend the Cox ring of $X$ (or its field of fractions)
  by several roots of homogeneous polynomials.
The number and the order of the roots is essentially bounded by the singularities of $Y$.

In the present article we extend this analogy to arbitrary \emph{Mori Dream Spaces}.
Suppose a normal algebraic variety $X$ over an algebraically closed field $\kk$ of characteristic $0$
   has a finitely generated divisor class group $\Cl(X)$.
Then we define its Cox ring to be:
\[
  S[X]:= \bigoplus_{[D] \in \Cl(X)} H^0(\ccO(D)).
\]
This $\kk$-vector space has a well defined ring structure (see \cite[\S 1.4]{arzhantsev_derenthal_hausen_laface_Cox_rings}),
  which makes it a multigraded ring with a grading by $\Cl(X)$.
We say $X$ is a \emph{Mori Dream Space} (or \emph{MDS}), if $S[X]$ is a finitely generated $\kk$-algebra.
In this situation there is a natural quotient map
$\pi_X \colon \Spec S[X] \dashrightarrow X$ by the action of the quasitorus $G_X= \Spec \kk[\Cl(X)] = \Hom(\Cl(X), \kk^*)$.

Affine and projective spaces and normal toric varieties are Mori Dream Spaces.
In these cases, the Cox ring is always a polynomial ring, but the grading vary.
In fact, the property that $S[X]$ is a polynomial ring characterises toric varieties,
  see \cite{kedzierski_wisniewski_Jaczewski_theorem_revisited}
  for a recent treatment of this characterisation.
Further examples of Mori Dream Spaces include projectivisations of rank two toric vector bundles~\cite{gonzalez_projectivised_rk_2_vb_are_MDS},
  del Pezzo surfaces~\cite{batyrev_popov_cox_ring_of_del_Pezzo_surf}
  and some K3 surfaces~\cite{artebani_hausen_laface_Cox_ring_of_K3}.
Every log Fano variety is a Mori Dream Space by~\cite{birkar_cascini_hacon_mckernan_minimal_model_for_log_general_type},
  in particular the moduli spaces of pointed rational curves $\overline{M}_{0,n}$
  have finitely generated Cox rings for $n\le 6$.
Minimal resolutions of surface quotient singularities are interesting non-projective examples of Mori Dream Spaces
  \cite{donten_bury_Cox_ring_of_minimal_resolutions_of_surface_quot_sings}.

In brief, our main result is the following.
\begin{thm}\label{thm_main_intro}
  Suppose $X$ and $Y$ are Mori Dream Spaces, and $\varphi\colon X \dashrightarrow Y$ is a rational map.
  Then there exists a description of $\varphi$ in terms of Cox coordinates, that is a multi-valued map
  \[
     \Phi \colon \Spec S[X] \multito \Spec S[Y]
  \]
  such that for all points $x \in X$ and  $\xi$ such that $\pi_X(\xi)=x$ and $\varphi$ is regular at $x$,
          the composition $\pi_Y(\Phi(\xi))$ is a single point $ \varphi(x) \in Y$.
\end{thm}

The precise definition of \emph{multi-valued map} is presented in Section~\ref{sec_multi-valued}.
It follows the convention of
   \cite[\S3.2]{brown_jabu_maps_of_toric_varieties}, where the same notion is used in the setup of toric varieties,
   that is the special case $\Spec S[X] =\kk^m$ and $\Spec S[Y] =\kk^n$.
The theorem is effective in the sense, that the proof shows how to construct the description.

Similar statement for regular maps between $\QQ$-factorial Mori Dream Spa\-ces
   has been independently obtained by Andreas Hochenegger and Elena Martinengo \cite{hochenegger_martinengo_maps_of_MDS}
Their approach is to use the language of Mori Dream stacks \cite{hochenegger_martinengo_MD_stacks}.
They use the technique of root constructions, which is parallel to the multi-valued maps in this article.

Theorem~\ref{thm_main_intro} is proved in several steps, concluding in Theorem~\ref{thm:existence_of_complete_descriptions} in Section~\ref{sec_complete_descriptions}.
It is a direct generalisation of the results of \cite{brown_jabu_maps_of_toric_varieties}.
Admittedly, a large part of this article adapts \cite[Sections 3 and 4]{brown_jabu_maps_of_toric_varieties} to the more general setting.
We heavily exploit the theory of Cox rings described in the book \cite{arzhantsev_derenthal_hausen_laface_Cox_rings} in order to make the generalisation possible.

%%%THINK: we may consider to extend the intro if needed.

\subsection*{Acknowledgements}

We would like to thank Andreas Hochenegger, Elena Martinengo and Jaros{\l}aw Wi{\'s}niewski for their hints, comments and corrections.
We thank IMPAN and the particiants of the minisemester miniPAGES 04-06.2016 for an excelent and inspiring atmosphere during the event.

\section{Preliminaries}

Fix an algebraically closed base field $\kk$ of characteristic $0$.
Suppose $X$ is a normal algebraic variety over $\kk$ with a finitely generated class group $\Cl(X)$.
We also assume that there are no non-constant global invertible functions on $X$, that is $H^0(X,\ccO_X^*)=\kk^*$.
Then $\ccS_X$, the \emph{Cox sheaf} of $X$, is the $\ccO_X$-algebra defined as
\[
  \ccS_X=\bigoplus_{[D]\in \Cl(X)} \ccO_X(D),
\]
and the \emph{Cox ring} of $X$ is the ring of global sections $S[X]=H^0(X,\ccS_X)$.
The construction of the multiplication structures in $\ccS_X$ and $S[X]$ is slightly delicate, especially in the case when $\Cl(X)$ is not torsion free,
  see~\cite[\S{}1.4.2]{arzhantsev_derenthal_hausen_laface_Cox_rings} for the details.

Assume $X$ is a \emph{Mori Dream Space}, that is $S[X]$ is finitely generated.
Then $S[X]$ is a domain by \cite[Thm~1.5.1.1]{arzhantsev_derenthal_hausen_laface_Cox_rings}.

We begin with listing the notation that will be used throughout the article.
\begin{itemize}
 \item The field of fractions of the Cox ring of $X$ will be denoted by $S(X)$.
       The contrast between the ring $S[X]$ and field $S(X)$ is analogous to the standard notation for polynomial ring
          $\kk[\fromto{x_1}{x_m}]$ and rational function field $\kk(\fromto{x_1}{x_m})$.
       Elements of $S[X]$ will be called \emph{regular sections} on $X$,
          elements of $S(X)$ will be called \emph{rational sections} on $X$.
 \item $K(X)$ is the field of rational functions on $X$.
       Elements of $K(X)$ will be called \emph{rational functions} on $X$.
 \item We fix an algebraic closure of the field $S(X)$ and we denote it by $\overline{S(X)}$.
 \item Since $S[X]$ is a graded ring, we can define the degree $0$ subfield of $S(X)$:
      \[
        S(X)^0=\set{\frac{f}{g}\in S(X)\;|\; f,g\in S[X], \text{ homogeneous, } \deg f=\deg g}.
      \]
      Note that $S(X)^0 = K(X)$ by Proposition~\ref{prop_degree_0_is_the_rational_function_field}.
  \item We say $q \in S(X)$ is homogeneous if it is a ratio of homogeneous elements of $S[X]$,
          i.e.~$q=\frac{f}{g}$ with $f, g$ homogeneous.
          Nevertheless note that $S(X)$ is not itself a graded ring.
  \item By $G_X$ we denote $\Spec\kk [\Cl(X)]$, the \emph{characteristic (quasi)torus}, whose character group is $\Cl(X)$.
  \item Let $\widehat{X} := \SheafySpec_X \ccS_X$, the relative spectrum of the sheaf of algebras $\ccS_X$.
        The variety $\widehat{X}$ is called the \emph{characteristic space} of $X$.
  \item Let $\overline{X}= \Spec S[X]$ be the \emph{total coordinate space} of $X$ (or a \emph{Cox cover}).
  \item Both $\widehat{X}$ and $\overline{X}$ admit the action of $G_X$
           determined by the grading of $\ccS_X$ and $S[X]$.
        The natural map $\widehat{X} \to \overline{X}$ is a $G_X$-equivariant open embedding \cite[Construction~1.6.3.1]{arzhantsev_derenthal_hausen_laface_Cox_rings}.
        The complement $\irrel{X} :=  \overline{X} \setminus \widehat{X}$ is called the \emph{irrelevant locus}.
  \item The \emph{irrelevant ideal} $B_X$ is the homogeneous ideal in $S[X]$ defining $\irrel{X} \subset \overline{X}$.
  \item Let $\pi_X\colon \overline{X} \dashrightarrow X$ be the natural rational map
           obtained by composing the inverse of the open embedding $\widehat{X} \hookrightarrow \overline{X}$
           and the relative spectrum map $\widehat{X} \to X$.
        Note that $\pi_X|_{\widehat{X}}$ is regular, surjective and it is the good quotient of $\widehat{X}$ by $G_X$ \cite[Construction~1.6.1.3]{arzhantsev_derenthal_hausen_laface_Cox_rings}.
  \item Suppose $A\subset X$ is a closed subscheme.
        Then the \emph{ideal} of $A$ is the homogeneous ideal $I(A) \subset S[X]$ generated by all sections vanishing on~$A$.
\end{itemize}

We will freely use the following equivalence:

\begin{prop}\label{prop_degree_0_is_the_rational_function_field}
   The rational function field $K(X)$ of a Mori Dream Space $X$ is naturally isomorphic to the field $S(X)^0$ of degree $0$ rational sections.
\end{prop}

\begin{prf}
   Suppose $f$ and $g$ are homogeneous elements of $S[X]$ of the same degree. Then $f, g \in H^0(\ccO_X(D))$ for a divisor $D$ on $X$.
   By definition
   \[
     H^0(\ccO_X(D)) = \set{q\in K(X) \mid D+(q) \ge 0}.
   \]
   Suppose $f$ corresponds to $q_1 \in K(X)$ and $g$ corresponds to $q_2 \in K(X)$ in this definition.
   Then
   \[
       S(X)^0 \ni \frac{f}{g} = \frac{q_1}{q_2} \in K(X),
   \]
   and this correspondence naturally defines a homomorphism $S(X)^0 \to K(X)$.
   It remains to verify that this map is surjective.

   Take a non-zero $q \in K(X)$. It defines a principal divisor $(q) \subset X$.
   Write $(q) = D^+ - D^-$, where both $D^+$ and $D^-$ are effective linearly equivalent divisors.
   Hence there exist homogeneous $f,g \in S[X]$ of the same degree with $D^+=(f)$ and $D^- = (g)$.
   Thus $f/g \in K(X)$ and $\frac{q}{f/g}$ is a rational function defining a trivial divisor, hence $q$ is a rescaling of $f/g$ by a globally invertible regular function,
     which must be constant by our assumptions.
\end{prf}

\subsection{Gradings and fields}

Let $M$ be a finitely generated abelian group, and $S = \bigoplus_{m \in M} S^m$ an $M$-graded domain.
Let $\FF$ be the field of fractions of $S$, so that we have $S \subset \FF$.
Consider the subfield of degree $0$ elements
\[
   \FF^0 = \set{\frac{f}{g} : f,g \in S,\text{ homogeneous, } \deg f - \deg g = 0}.
\]
If $A\subset \FF$ is a subset, then by $\FF^0(A)$ we denote the smallest subfield of $\FF$ containing $\FF^0$ and $A$, i.e.~the field generated over $\FF^0$ by $A$.

Given a subgroup $M'$ of $M$ we may also consider a coarser grading on $S$, namely by the finitely generated group $M/M'$:
If $[m] \in M/M'$, then $S_{[m]}=\bigoplus_{m'\in M'} S_{m+m'}$ with respect to the coarser grading.
We define the subfield of degree zero elements with respect to the coarser grading:
\[
  \FF^{M'} = \set{\frac{f}{g} : f,g \in S_{[m]} \text{ for some coarse degree $[m] \in M/M'$}}
\]
(we stress that the coarse degrees of $f$ and $g$ are the same in the formula above).
\begin{lemma}\label{lem_extending_fields_by_homogeneous_elements}
   Suppose $A \subset S$ is a subset of the ring consisting only of homogeneous elements (in particular, $A$ may be infinite).
   Let $\langle\deg{A}\rangle \subset M$ be the subgroup of $M$ generated by the degrees of elements of $A$. (For consistence, we assume $\deg 0 = 0$.)
   Then
   \[
       \FF^0(A) = \FF^{\langle\deg{A}\rangle}.
   \]
   In particular, there is a finite subset $A'\subset A$ such that $\FF^0(A) = \FF^0(A')$.
\end{lemma}
\begin{prf}
   We always have $\FF^0 \subset \FF^{\langle\deg A\rangle}$ and $A\subset \FF^{\langle\deg A \rangle}$, so:
   \begin{equation}\label{equ_easy_field_inclusion}
      \FF^0(A) \subset \FF^{\langle\deg A \rangle}.
   \end{equation}
   To show the reverse inclusion $\FF^0(A) \supset \FF^{\langle\deg A \rangle}$, we first suppose
      the set $A$ consists of a single homogeneous element $A = \set{a}$.
   If $a=0$, there is nothing to prove.
   Otherwise, pick $q\in \FF^{\langle\deg A \rangle}$, i.e.
   \[
     q= \frac{f_i + f_{i+1}+ \dotsb + f_{j}}{g_k + g_{k+1}+ \dotsb + g_{l}}
   \]
   where $i,j,k,l \in \ZZ$ and $f_n, g_n \in  S_{n \deg a}$.
   In particular, $\frac{f_n}{a^n} \in \FF^0$ and analogously $\frac{g_n}{a^n} \in \FF^0$.
   Thus:
   \[
     q= \frac{\frac{f_i}{a^i} a^i + \dotsb + \frac{f_j}{a^j} a^j}{\frac{g_k}{a^k} a^k + \dotsb + \frac{g_l}{a^l} a^l}
   \]
   which expresses $q$ in terms of $a$ and $\FF^0$. Thus $\FF^{\langle \deg a \rangle} \subset \FF^0(\set{a})$.

   Now assume $A$ is finite. We can show the statement inductively, by adding the generators one by one and applying the initial case (with $A$ equal to a single element).
   More precisely, let $A= A'\cup \set{a}$.
   By the inductive assumption $\FF^0(A')  = \FF^{\langle\deg(A')\rangle}$ and we have
   \[
     \FF^0(A)  = \left(\FF^0(A')\right)(\set{a}) = \FF^{\langle\deg(A')\rangle}(\set{a}).
   \]
   Considering the coarser grading by $M/\langle\deg(A')\rangle$ the field $\FF^{\langle\deg(A')\rangle}$ is really a subfield of degree $[0]$ elements
      (with respect to the coarser grading).
   Thus we apply the initial case to obtain $\FF^{\langle\deg(A')\rangle}(\set{a}) = \FF^{\langle\deg(A)\rangle}$.

   Finally, suppose $A$ is arbitrary. Since $M$ is finitely generated and abelian, also $\langle\deg A \rangle$,
      which is a subgroup of $M$, is finitely generated.
   Pick a finite subset $A' \subset A$, such that $\langle\deg A' \rangle = \langle\deg A \rangle$. Using this equality and the case when $A$ is finite we obtain
   \[
      \FF^0(A) \stackrel{\text{\eqref{equ_easy_field_inclusion}}}{\subset} \FF^{\langle\deg A \rangle} = \FF^{\langle\deg A' \rangle}
               = \FF^{0}(A') \stackrel{A'\subset A}{\subset} \FF^0(A).
   \]
   Thus all the inclusions above are equalities, i.e.~$\FF^0 (A) = \FF^{\deg(A)} = \FF^{0}(A')$, concluding the proof.
\end{prf}

\subsection{Maps of fields with kernel}\label{sec_maps_of_fields_with_kernel}

Given two algebraic varieties $X$ and $Y$ over an algebraically closed field $\kk$ we consider a rational map $\varphi\colon X \dashrightarrow Y$.
If $\varphi$ is dominant, then such map is uniquely determined by a homomorphism of fields $\varphi^*\colon K(Y) \to K(X)$.
In this section we are interested in the algebraic interpretation of the case when $\varphi\colon X \dashrightarrow Y$ is not necessarily dominant, compare to \cite[Prop.~2.14(i),(ii)]{brown_jabu_maps_of_toric_varieties}.

One approach is to replace $Y$ with the closure of the image $Z:=\overline{\varphi(X)}$.
Then $\varphi\colon X \dashrightarrow Z$ is a dominant and there is a homomorphism $\varphi^*\colon K(Z) \to K(X)$.
However, this does not help very much to describe, for example, closed embeddings.
Here we are seeking a uniform description in terms of $K(Y)$, as we are going to extend this description to Cox ring $S[Y]$.
For this purpose we define the notion of a \emph{map of fields with kernel}.

\begin{defin}\label{def_map_of_fields_with_kernel}
   Suppose $\FF$ and $\GG$ are fields. A \emph{map of fields with kernel} (denoted $\Phi^* \colon \FF \dashrightarrow \GG$)
     is a ring homomorphism $\Phi^* \colon R \to \GG$, where:
     \begin{itemize}
        \item $R \subset \FF$ is a subring which generates $\FF$ (i.e.~$\FF$ is the smallest subfield of $\FF$ containing $R$),
        \item if $f\in R$ and $f^{-1} \in \FF \setminus R$, then $\Phi^*(f) =0$.
     \end{itemize}
\end{defin}

Thus if $\varphi\colon X \dashrightarrow Y$ is a rational map, then let $R\subset K(Y)$ be the local ring of the scheme-theoretic point,
   whose closure is $\overline{\varphi(X)}$, or equivalently, the set of those rational functions on $Y$, whose set of poles does not contain $\varphi(X)$.
Clearly $R$ generates $\FF$, and $\varphi^* \colon R \to K(X)$ is a well defined homomorphism.
If $f \in R$ and $f^{-1} \in \FF \setminus R$, then $\varphi(X)$ is contained in the set of poles of $f^{-1}$ and but is not contained in the set of poles of $f$,
   hence $f$ must be $0$ along $\overline{\varphi(X)}$.
Thus $\varphi^* \colon K(Y) \dashrightarrow K(X)$ is a map of fields with kernel.

% We can work with maps of fields with kernel in a similar way as with the homomorphisms of fields.
\begin{lemma}\label{lem_map_of_fields_with_kernel_can_be_defined_on_subrings}
   Let $\FF$, $\GG$ be fields, $R' \subset \FF$ be a subring generating $\FF$, and let $\Phi^*\colon R' \to \GG$ be a ring homomorphism.
   Then there is a unique map of fields with kernel $\FF \dashrightarrow \GG$ extending $\Phi^*$.
\end{lemma}

\begin{prf}
   Consider the kernel of $\Phi^* \colon R' \to \GG$ and the localisation $R:=(R')_{\ker_{\Phi^*}}$ of $R'$ in this ideal.
   Then $\Phi^*$ extends uniquely to $R$ and $R$ generates $\FF$.
   If $f\in R$ and $f^{-1} \notin R$, then $f \in \ker \Phi^*$, hence $\Phi^*f=0$ and $\Phi^*\colon \FF \dashrightarrow \GG$ is a map of fields with kernel.
\end{prf}

% \begin{lemma}\label{lem_maps_of_rings_with_kernel_can_be_extended}
%    Suppose $\FF$, $\FF^0$ and $\GG$ are fields, $\FF^0\subset \FF$ and $\Phi^*\colon \FF^0 \dashrightarrow \GG$ is a map of fields with kernel,
%       and $R^0 \subset \FF^0$ is the ring where $\Phi^*$ is defined as a ring homomorphism.
%    Assume $\FF = \FF^0(a)$ for a single element $a\in \FF$.
%    In addition suppose $b\in \GG$ and:
%    \begin{itemize}
%     \item  If $a$ is transcendental element over $\FF^0$ (or equivalently $\FF \simeq \FF^0(t)$),
%             then $b\in \GG$ is any element.
%     \item  If $a$ satisfies a minimal polynomial relation $F(a)=0$,
%                  where $F$ is a monic polynomial in one variable with coefficients in $R^0$,
%              then $b$ is an element in $\GG$ satisfying $\Phi^*F(b)=0$.
%    \end{itemize}
%    Then we can define $\Psi^*\colon \FF \dashrightarrow \GG$ such that $\Psi^*|_{\FF^0} = \Phi^*$ and $\Psi^*(a)=b$.
% \end{lemma}
%
% \begin{prf}
%    The conditions on $a$ and $b$ allow us to extend the ring homomorphism $\Phi^*$ to $R^0[a] \to \GG$
%       (where $R^0[a]$ is the subring of $\FF$ generated by $R^0$ and $a$).
%    Thus the claim follows from Lemma~\ref{lem_map_of_fields_with_kernel_can_be_defined_on_subrings}.
% \end{prf}
%

\section{Multivalued sections and maps}\label{sec_multi-valued}

Suppose $X$ and $Y$ are Mori Dream Spaces, and consider their total coordinate spaces
   $\overline{X}= \Spec S[X]$ and $\overline{Y} = \Spec S[Y]$.

\begin{defin}
   A \emph{homogeneous multi-valued section} on $X$ is an element $\gamma\in\overline{S(X)}$ of the algebraic closure of
      the fraction field $S(X)$, such that $\gamma^r \in S(X)$ is a homogeneous rational section on $X$ for some positive integer $r$.
\end{defin}

We always pick minimal $r$, such that $\gamma^r = q\in S(X)$, and then we write $\gamma = \sqrt[r]{q}$.
Then a value of $\gamma$ at a point $\xi \in \overline{X}$ is $\gamma(\xi) = \set{\sqrt[r]{q(\xi)}} \subset \kk$,
   i.e.~the set of all $r$-th roots of $q(\xi)$ (assuming $\xi$ is not a pole of $q$).

\begin{prop}\label{prop:single_valued_sections}
  A homogeneous multi-valued section $\gamma \in \overline{S(X)}$ is in $S(X)$
  if and only if $\gamma(\xi)$ has exactly one element for a general $\xi \in \overline{X}$.
\end{prop}
\noprf

\begin{prop}\label{prop:image_has_constant_degree_and_irrational_part}
If $V$ is a $\kk$-vector space and $i \colon V \to \overline{S(X)}$
is a $\kk$-linear map whose image consists of only homogeneous multi-valued
sections, then there exists a homogeneous multi-valued section $\gamma \in
\overline{S(X)}$ and a $\kk$-linear map $j \colon V \to
S(X)$ whose image  consists of homogeneous elements of a constant
degree, and $ i(v) = j(v) \cdot \gamma$  for all  $v\in V$.
\end{prop}

See \cite[Cor.~2.20 and Prop.~3.6]{brown_jabu_maps_of_toric_varieties};
the main ingredient of the proof is \cite[Cor.~2.20]{brown_jabu_maps_of_toric_varieties},
which is valid over any base field, the only assumption is that $\FF=S(X)$ contains all roots of unity.
On the other hand, the statement of \cite[Prop.~3.6]{brown_jabu_maps_of_toric_varieties}
  is analogous to Proposition~\ref{prop:image_has_constant_degree_and_irrational_part},
  but restricted to toric varieties (i.e.~$S(X)$ is a field of rational functions).
However this assumption is never exploited in the simple argument.

A regular map between affine varieties is merely a geometric interpretation
   of a homomorphism of respective coordinate rings.
In the same spirit we want to give a geometric interpretation to a more general homomorphism,
which is composed of homogeneous multi-valued sections.

For this purpose, fix $\fromto{y_1}{y_N}$, a finite number of homogeneous generators of the Cox ring $S[Y]$.
We will see later, in Section~\ref{sec_homogeneity_condition},
  that the particular choice of the generators does not matter from the point of view of our needs.
\begin{defin}
   A \emph{multi-valued map} $\Phi \colon \overline{X}\multito \overline{Y}$ is a geometric interpretation of a ring homomorphism
     $\Phi^* \colon S[Y] \to  \overline{S(X)}$, which maps the generators $y_i$ to homogeneous multi-valued sections of $X$.
   That is, for all $i$, there exists $r_i$, such that $\Phi^*(y_i^{r_i})\in S(X)$ is a homogeneous rational section on $X$.
\end{defin}

The purely algebraic definition of $\Phi \colon \overline{X}\multito \overline{Y}$ has now immediate geometric consequences.
First, we define  $\Reg \Phi \subset \overline{X}$ the \emph{regular locus} of $\Phi$, to be the open dense subset,
    complementary to the zero sets of denominators appearing in $(\Phi^*y_i)^{r_i}$.
That is, we write $\Phi^*y_i = \sqrt[r_i]{\frac{f_i}{g_i}}$ for homogeneous regular sections $f_i, g_i \in S[X]$,
and we assume $f_i$ and $g_i$ have no common homogeneous factors.
This is possible by the graded factoriality of the Cox ring, see~\cite[Theorem~1.5.3.7]{arzhantsev_derenthal_hausen_laface_Cox_rings}.
Then:
\[
   \Reg \Phi := \set{\xi \in \overline{X} \mid \forall_i \ g_i(\xi) \ne 0} = \overline{X} \setminus Z(g_1 \dotsm g_N).
\]
\begin{cor}\label{cor_regularity_locus_is_affine}
   $\Reg \Phi$ is the complement of a $G_X$-invariant divisor in $\overline{X}$ and it is affine.
\end{cor}

Suppose $\Phi \colon \overline{X}\multito \overline{Y}$ is a multi-valued map.
Since $S[Y]$ generates $S(Y)$, the homorphism $\Phi^*$
  determines a unique map of fields with kernel $\Phi^* \colon S(Y) \dashrightarrow \overline{S(X)}$
  by Lemma~\ref{lem_map_of_fields_with_kernel_can_be_defined_on_subrings}.
Specifically, pick $q = \frac{f}{g} \in S(Y)$  with $f,g \in S[Y]$.
Assume $\Phi^*g \ne 0$, then $\Phi^* q:= \frac{\Phi^*f}{\Phi^*g}$ is well defined,
   and thus $\Phi^*$ naturally extends to:
\begin{equation} \label{equ_extending_Phi_to_rational_functions}
   \Phi^* \colon R \to \overline{S(X)}, \text{ where $R$ is the subring }
   R = \set{q =\tfrac{f}{g} \in S(Y) \mid \Phi^*g \ne 0}.
 \end{equation}

\begin{rem}\label{rem_subring_R0}
   In particular, $R$ contains a subring $R^0 \subset R$ of quotients of degree $0$, hence $R^0\subset K(Y)$.
\end{rem}

Let $\Gamma(\Phi)$ be the subring of $\overline{S(X)}$ generated by $S[X][g^{-1}]$, and $\Phi^*(S[Y])$, where $g=g_1\dotsm g_N$,
  so that $\Reg \Phi = \Spec S[X][g^{-1}]$.
Then we have the following algebraic morphisms:
\[
 \overline X  \supset  \Reg \Phi \stackrel{p_{\Phi}}{\longleftarrow} \Spec \Gamma(\Phi) \stackrel{q_{\Phi}}{\longrightarrow} \Spec S[Y].
\]
The first map $p_{\Phi}$ is the morphism of affine varieties corresponding to the inclusion
   $p_{\Phi}^* \colon S[X][g^{-1}] \to \Gamma(\Phi)$.
The second map $q_{\Phi}$ corresponds to $q_{\Phi}^*  = \Phi^*  \colon S[Y] \to \Gamma(\Phi)$.
Using these two maps we define the (set-theoretic) image and preimage under a multi-valued map:
\begin{defin}\label{def_image_and_preimage}
   Suppose $Z \subset \overline{X}$ and $W \subset \overline{Y}$ are two subsets. Then define:
   \[
      \Phi(Z):= q_{\Phi}(p_{\Phi}^{-1}(Z)) \quad \text{and} \quad \Phi^{-1}(W):=  p_{\Phi}(q_{\Phi}^{-1} (W)).
   \]
\end{defin}
In particular, the following are immediate consequences of the definition:
\begin{prop}
In the setting of Definition~\ref{def_image_and_preimage}:
   \begin{enumerate}
      \item $\Phi(Z) = \Phi(Z \cap \Reg \Phi)$
      \item $\Phi^{-1}(W)\subset \Reg \Phi$.
      \item $\Phi(Z) = \emptyset$ if and only if $Z \subset \overline{X} \setminus \Reg \Phi$,
            i.e.~$\Reg \Phi$ is the set of points, where the image under $\Phi$ is defined.
      \item For $g\in S[Y]$ we have $\Phi^*g = 0$ if and only if $\Phi(\overline{X}) \subset (g=0)$.
      \item Let $I = \ker \Phi^*$.
            Then the ring $R$ from Equation~\eqref{equ_extending_Phi_to_rational_functions}
            is equal to the (non-homogeneous) localisation $S[Y]_{I}$.
   \end{enumerate}
\end{prop}

Observe that $p_{\Phi} \colon\Spec \Gamma(\Phi)\to \Reg \Phi$ is a finite morphism. In particular, $p_{\Phi}$ is closed.
Thus (since $q_{\Phi}$ is continuous):
\begin{prop}\label{prop_preimage_of_closed_is_closed}
   For a closed subset $W \subset \overline{Y}$ the preimage $\Phi^{-1}(W)$ is closed in $\Reg \Phi$.
\end{prop}

\begin{rem}
   It is also possible to show that the preimage (under a  multi-valued map $\Phi$) of an open subset is open, thus $\Phi$ has very much of properties of continuous mappings.
   This statement requires a more refined definition of $\Gamma_{\Phi}$, and we postpone it until the second part of these series of papers.
   See also \cite[\S3.3, \S3.4, Prop.~3.17]{brown_jabu_maps_of_toric_varieties}.
\end{rem}

\section{Homogeneity conditions}\label{sec_homogeneity_condition}

Analogously to \cite[Prop.~4.6]{brown_jabu_maps_of_toric_varieties},
   we prove the equivalence of several	 conditions, jointly referred to as \emph{homogeneity conditions}.

Recall the setup: $X$ and $Y$ are Mori Dream Spaces, with total coordinate spaces $\overline{X} =\Spec S[X]$ and $\overline{Y} = \Spec S[Y]$.
We pick some homogeneous generators $\fromto{y_1}{y_N}$ of $S[Y]$.

\begin{prop} \label{prop:equivalence_of_homog_conds}
   Suppose $\Phi \colon \overline{X} \multito \overline{Y}$ is a multi-valued map and
     consider the set
   \[
      L=\set{ y_i \mid i \in \setfromto{1}{N} \text{ and }\ \Phi^* y_i \ne 0}
   \]
   of the fixed generators the Cox ring $S[Y]$ that pull back nontrivially under $\Phi$.
   Further let $\FF \subset S(Y)$ be the subfield generated by $L$,
     and let $\FF^0 := \FF \cap S(Y)^0 \subset K(Y)$.
   The following conditions are equivalent:
   \begin{enumerate}
      \renewcommand{\theenumi}{\textnormal{(A\arabic{enumi})}}
      \item
         \label{item:homogeneity_cond_all}
         If $q\in S(Y)$ is homogeneous and $\Phi^*q$ is defined (i.e.~the pullback of the denominator is non-zero), then
         $\Phi^*q$ is a homogeneous multi-valued section on $X$.
      \item
         \label{item:homogeneity_cond_on degree_0}
         If $q\in K(Y)$ and $\Phi^*q$ is defined (i.e.~$q \in R^0$ in the sense of Remark~\ref{rem_subring_R0}),
           then $\Phi^*q\in K(X)$.
      \item
         \label{item:homogeneity_cond_monomial}
         There exist $l_1,\dotsc,l_k$, elements of the field $\FF^0$, generating it as a field extension of $\kk$
            such that $\Phi^* l_i$ are elements of $S(X)^0 = K(X)$.
      \item
         \label{item:homogeneity_cond_points}
         $\Phi$ maps $G_X$-orbits into $G_Y$-orbits.
         More precisely, for all $\xi,\xi' \in \Reg\Phi$ with $\xi' \in G_X \cdot \xi$,
           if $\eta \in \Phi(\xi)$ and $\eta' \in \Phi(\xi')$ then $\eta' \in G_Y \cdot \eta$.
      \renewcommand{\theenumi}{\textnormal{(A\arabic{enumi}${}^\prime$)}}
      \addtocounter{enumi}{-1}
      \item
         \label{item:homogeneity_cond_points_on_open}
         $\Phi$ maps $G_X$-orbits meeting an open dense subset of $X$ into $G_Y$-orbits.
         More precisely,
           there exists an open dense subset $U \subset \Reg \Phi$, such that
           for all $\xi,\xi' \in U$ with $\xi' \in G_X \cdot \xi$,
           if $\eta \in \Phi(\xi)$ and $\eta' \in \Phi(\xi')$ then $\eta' \in G_Y \cdot \eta$.
   \end{enumerate}
\end{prop}

\begin{prf}
Suppose \ref{item:homogeneity_cond_all} holds for $\Phi$.
Let $R^0 \subset S(Y)$ be the $\kk$-vector subspace of
  homogeneous sections of degree $0$ for which the pull-back by $\Phi$ is defined as in Remark~\ref{rem_subring_R0}.
Denote the restriction of $\Phi^*$ to $R^0$ by
$i\colon R^0 \to \overline{S(X)}$. Since $i(1)=1$ is rational
and has degree $0$,
Proposition~\ref{prop:image_has_constant_degree_and_irrational_part} implies
that all elements of $i(R^0)$ are rational and of degree $0$.
Therefore \ref{item:homogeneity_cond_on degree_0} holds for $\Phi$.

Suppose \ref{item:homogeneity_cond_on degree_0} holds. Since
$\FF^0 \subset S(Y)^0$ (in fact $\FF^0 \subset R^0$),
any generating set $\fromto{l_1}{l_k}$ of $\FF^0$ satisfies
\ref{item:homogeneity_cond_monomial}.

Suppose \ref{item:homogeneity_cond_monomial} holds for $\Phi$; we show
that \ref{item:homogeneity_cond_all} holds.
Let $q\in S(Y)$ be any homogeneous rational section. Write
\[
   q = \frac{\mu_1 + \dotsb + \mu_{\alpha}}{\nu_1 + \dotsb + \nu_{\beta}},
\]
where the $\mu_i, \nu_j\in S[Y]$ are monomial expressions in the generators $y_{k}$
(i.e.~simply products of powers of the generators $\fromto{y_1}{y_N}$ with non-negative exponents)
with $\deg \mu_i =d_1$ and  $\deg \nu_j =d_2$ for all $i$ and~$j$. Assume that
$\Phi^*(\nu_1 + \dotsb + \nu_{\beta}) \ne 0$, so $\Phi^*q$ is defined.

We have to show $\Phi^*q$ is a homogeneous multi-valued section.
If the expression of $\mu_i$ as the product of the generators
  contains $y_k$, such that $\Phi^*y_k=0$, then $\Phi^*\mu_i =0$.
Let $q' = \frac{\mu_1 + \dotsb + \mu_{i-1} +\mu_{i+1} +\dotsb + \mu_{\alpha}}{\nu_1 + \dotsb + \nu_{\beta}}$.
Then $\Phi^*(q')= \Phi^*(q)$, hence we can replace $q$ with $q'$.
Repeating the operation if necessary (also for the denominator),
   we may assume that $\mu_i$ and $\nu_j$ are all monomial expressions in the generators from $L$ only, i.e.~$q \in \FF$.

Since $\mu_{i_1}/\mu_{i_2}$ is homogeneous of degree $0$, thus $\mu_{i_1}/\mu_{i_2} \in \FF^0$,
   i.e.~it is expressible in terms of the generators $l_i$.
Thus, using the assumption of \ref{item:homogeneity_cond_monomial},
the pull-back $\Phi^*(\mu_{i_1}/\mu_{i_2})$ is a non-zero homogeneous degree~$0$ rational section in $S(X)$.

Each $\Phi^*\mu_i$ is a homogeneous multi-valued section (it is a product of homogeneous multi-valued sections).
In particular, for every $i$,
\[
    \Phi^*(\mu_i) = f_i \cdot \gamma
\]
where $\gamma \in \overline{S(X)}$ is a fixed homogeneous multi-valued
section (independent of $i$) and $f_i \in K(X)$. So
\[
    \Phi^*(\mu_1  + \dotsb + \mu_{\alpha}) = (f_1 + \dotsb + f_{\alpha}) \gamma.
\]
Similarly,
$\Phi^*(\nu_1  + \dotsb + \nu_{\beta}) = (g_1 + \dotsb + g_{\beta}) \delta \not= 0$,
for some $\delta \in \overline{S(X)}$ and $g_j\in K(X)$.
So
\[
   \Phi^*(q) = h \cdot {\varepsilon}
\]
where $\varepsilon = \gamma/\delta \in \overline{S(X)}$ is homogeneous
and $h = (\sum f_i) / (\sum g_j) \in K(X)$.
So $\Phi^*(q)$ is a homogeneous multi-valued section and \ref{item:homogeneity_cond_all} holds.

It remains to prove the implications
\ref{item:homogeneity_cond_all}
$\Longrightarrow$ \ref{item:homogeneity_cond_points}
$\Longrightarrow$  \ref{item:homogeneity_cond_points_on_open}
$\Longrightarrow$ \ref{item:homogeneity_cond_on degree_0}.

    Suppose \ref{item:homogeneity_cond_all} holds.
    Let $\xi \in \Reg \Phi$  and consider $G_X \cdot \xi$.
    The claim of \ref{item:homogeneity_cond_points} is that $\Phi(G_X \cdot\xi)$ is contained in one $G_Y$ orbit.
    Let $I \subset S[Y]$ be the homogeneous ideal generated by:
    \[
       I := \langle f  \in S[Y] \;|\; f \text{ is homogeneous and } (\Phi^*f)(\xi) = 0   \rangle
    \]
    If $f \in I$, then $\Phi^*f (\xi') =0$ for all $\xi'$ in the orbit $G_X \cdot \xi$.
    Thus $\Phi(G_X \cdot \xi)$ is contained in the ($G_Y$-invariant) zero locus $T:=Z(I)\subset \overline{Y}$.

    We claim that for any $\eta' \in \Phi(G_X \cdot \xi)$ the orbit $G_Y\cdot \eta'$ is dense in $T$.
    Suppose on the contrary, that $\overline{G_Y\cdot \eta'} \subsetneqq T$,
       and thus there exists a homogeneous $f \in S[Y]$, which is not in $I$,
       such that $f(\eta')=0$.
    Let $\xi' \in G_X \cdot \xi$ be such that $\eta'\in \Phi(\xi')$.
    Then $\Phi^*f (\xi') = f(\eta')=0$. But $\Phi^*f$ is homogeneous, hence also $\Phi^*f (\xi) =0$,
       a contradiction, since we assumed $f \notin I$.

    A dense orbit, if exists, is unique (even if $T$ is not irreducible!).
    Hence by above, the action of $G_Y$ on $T$ has a dense orbit and
      $\Phi(G_X \cdot \xi)$ is contained in this orbit.
    This completes the proof of  \ref{item:homogeneity_cond_points}.

    If \ref{item:homogeneity_cond_points} holds, then clearly \ref{item:homogeneity_cond_points_on_open} holds too.

    Finally, suppose \ref{item:homogeneity_cond_points_on_open} holds.
    Let $q\in K(Y)$ be such that $\Phi^*q$ is defined.
    Suppose $\xi \in U$ is general.
    The possible values taken by $\Phi^*q$ at $\xi$ are simply those values taken by $q$
       at the points of the image set $\Phi(\xi)$.
    Setting $\xi'=\xi$ in \ref{item:homogeneity_cond_points_on_open}
       shows that $\Phi(\xi)$ is contained in a single $G_Y$-orbit,
    and so
    $
       \Phi^*q(\xi) = \set{q(\eta) \mid \eta \in \Phi(\xi)}
    $
    is a single number.
    Therefore $\Phi^* q \in S(X)$ by Proposition~\ref{prop:single_valued_sections}.
    For arbitrary $\xi$, $\xi'$ as in \ref{item:homogeneity_cond_points_on_open},
    \[
       (\Phi^*q)(\xi) = q(\Phi(\xi)) = q(\Phi(\xi')) = (\Phi^*q)(\xi')
    \]
    since $q$ is constant on $G_Y$-orbits. That is, $\Phi^* q$ is constant on a general $G_X$ orbit,
       so $\Phi^* q$ is $G_X$-invariant. This proves \ref{item:homogeneity_cond_on degree_0}.
\end{prf}

Note that by the intrinsic nature of Condition~\ref{item:homogeneity_cond_all}
   also all other Conditions~\ref{item:homogeneity_cond_on degree_0}, \ref{item:homogeneity_cond_monomial}, \ref{item:homogeneity_cond_points}, \ref{item:homogeneity_cond_points_on_open} are independent of the choice of the generators $y_i$ of $S[Y]$.

\section{Relevance condition}\label{sect_relevance}

We define and briefly analyse the relevance condition, analogous to \cite[Prop.~4.7]{brown_jabu_maps_of_toric_varieties}.

\begin{defin}
   Let $\Phi \colon \overline{X} \multito \overline{Y}$ be a multi-valued map.
   We say that $\Phi$ satisfies the \emph{relevance condition}, or $\Phi$ is \emph{relevant}, if
   \begin{enumerate}
      \renewcommand{\theenumi}{\textnormal{(B1)}}
      \item \label{item:rel_cond_points}
          the image of $\Phi$ is not contained in the irrelevant locus of $Y$,
   \end{enumerate}
   or equivalently,
   \begin{enumerate}
      \renewcommand{\theenumi}{\textnormal{(B2)}}
      \item \label{item:rel_cond_kernel}
          $\ker \Phi^*$ does not contain the irrelevant ideal $B_Y$ of $Y$.
   \end{enumerate}
\end{defin}

An important class of open subsets in $Y$ arise as complements of divisors.
For some of the arguments below we are interested in the the open affine subsets of $Y$. 
It is possible to show, using algebraic version of Hartog's theorem, that every open affine subset is a complement of a divisor 
  of the form  $Y \setminus \Supp (f)$ for a homogeneous regular section $f\in S[Y]$.
Here we show a rather weak version of this statement that we will use.
For a homogeneous regular section $f \in S[Y]$ we define 
   $\overline{U}_f = \overline{Y} \setminus \set{f= 0}$,
     i.e.~$\overline{U}_f = \Spec S[Y][\frac{1}{f}]$ is an open affine $G_Y$-invariant subset of $\overline{Y} = \Spec S[Y]$.
   Moreover, $\overline{U}_f$ is contained in the characteristic space $\widehat{Y} = \overline{Y} \setminus Z(B_Y)$.

\begin{lemma}\label{lem_open_affines_U_f}
   Suppose $f \in B_Y \subset S[Y]$ is a homogeneous regular section in the irrelevant ideal. 
   Then $U_f = \pi_Y(\overline{U}_f)$ is an open affine subset of $Y$.
   Moreover, $U_f = \Spec S[Y]_{(f)}$, 
     the homogeneous localisation of $S[Y]$ in $f$, 
     that is the ring of degree $0$ quotients $\frac{g}{f^d}$ for some homogeneous $g \in S[Y]$,
     and nonnegative integer $d$ such that $\deg g = d\cdot \deg f$.

   In addition, open affine subsets of the form $U_f$ (for homogeneous $f\in B_Y$) cover $Y$.
\end{lemma}
\begin{prf}
   The geometric invariant theory (GIT) guarantees that  $U_f$ is an open affine subset of $Y$.
   Indeed, it is open by \cite[Cor.~3.7]{swiecicka_good_quotients_survey_wykno}.
   It is affine since $\overline{U}_f \to U_f$ is a good quotient and thus
     \[
       U_f = \Spec \ccO_{\overline Y}(\overline{U}_f)^{G_Y} = \Spec S[Y]_{(f)}.
     \]
     
   To see that such sets cover $Y$, take $y = \pi_Y(\eta) \in Y$ for some $\eta \in \widehat{Y}$.
   Since $B_Y$ is a homogeneous ideal and $\eta \notin Z(B_Y) \subset \overline{Y}$,
     there is a homogeneous element $f$ of $B_Y$ which does not vanish at $\eta$.
   Then $\eta \in \overline{U}_f$ and $y =\pi(\eta) \in \pi(\overline{U}_f) = U_f$.
\end{prf}

In the setting of the relevance condition we have the following conclusions.

\begin{lemma}\label{lem_relevant_then_R0_generates_field}
  In the setting as above, let the ring $R^0$ be as in Remark~\ref{rem_subring_R0}.
  \begin{itemize}
    \item If $\Phi$ is relevant, then $R^0$ generates $K(Y)$ as a field,
             or, in other words, the map of fields with kernel $\Phi^*\colon S(Y) \dashrightarrow \overline{S(X)}$
             restricts to a map of fields with kernel $\Phi^*\colon K(Y) \dashrightarrow \overline{S(X)}$.
    \item If $\Phi$ satisfies both homogeneity and relevance conditions,
             then $\Phi^*$ restricts to a map of fields with kernel $\Phi^*\colon K(Y) \dashrightarrow K(X)$.

  \end{itemize}
\end{lemma}

\begin{prf}
  Take a homogeneous regular section $f\in B_Y$ such that $f\notin \ker \Phi^*$.
  Consider $U_f = \Spec S[Y]_{(f)} \subset Y$ as in Lemma~\ref{lem_open_affines_U_f}.
  The first part of the lemma follows from the fact that $R^0=(S[Y]_{\ker \Phi^*})^0$ 
  and that the latter contains $S[Y]_{(f)}$ which generates the field $K(Y)$.
  The equivalent phrasing in terms of maps of fields with kernel follows immediately from Lemma~\ref{lem_map_of_fields_with_kernel_can_be_defined_on_subrings}.
  The second item follows from the homogeneity condition \ref{item:homogeneity_cond_on degree_0}.
\end{prf}

\section{Description of maps}
\label{sect:descriptions}

In this section we exploit the language introduced so far to explain when
  a multi-valued map $ \overline{X} \multito \overline{Y}$ describes a rational map $X\dashrightarrow Y$.
We characterise those multi-valued maps that describe some rational map in terms of the homogeneity and relevance conditions.
Finally, we prove that every rational map between Mori Dream Spaces has such a description.

Let $\Phi\colon \overline{X} \multito \overline{Y}$ be a multi-valued map.
Define a subset $\Reg_Y\Phi \subset \Reg \Phi$,
the locus where $\pi_Y \circ \Phi$ is a well-defined map of sets:
\[
  \Reg_Y\Phi:=
    \set{ \xi \in \Reg \Phi \mid \# \pi_Y(\Phi(\xi)) = 1}.
\]
This locus $\Reg_Y\Phi$ may be empty.
On the other hand, if $\Reg_Y\Phi$ contains a nonempty open subset,
   then we regard $\Phi$ as being adapted to $Y$;
under this assumption, it makes sense to ask where $\Phi$ \emph{agrees} with a rational map $X\dashrightarrow Y$.

\begin{defin}
 Given a multi-valued map $\Phi\colon \overline{X} \multito \overline{Y}$
 and a rational map $\varphi\colon X \dashrightarrow Y$, in the notation above,
 the \emph{agreement locus of $\Phi$ and $\varphi$} is
 \[
   \agree{\Phi}{\varphi} = \set{
      \xi\in \Reg_Y\Phi \cap {\pi_X}^{-1}(\Reg\varphi) \mid
      \pi_Y\circ\Phi(\xi) = \varphi\circ\pi_X(\xi) }.
 \]
\end{defin}
In other words, the agreement locus is the set of points where both
compositions $\pi_Y \circ \Phi$ and $\varphi \circ \pi_X$ are well-defined
as maps of sets and they have the same values.

\begin{defin}\label{defin:description}
	We say \emph{$\Phi$ is a description of $\varphi$},
	or that \emph{$\Phi$ describes $\varphi$},
	if $\agree \Phi \varphi$ contains an open dense subset of $\overline{X}$.
\end{defin}

Note that the definition of a description is a point-wise definition.
However, it has strong algebraic consequences.
Firstly, a description uniquely determines the map $\varphi$ by \cite[Lem.~I.4.1 and Def. on page 24]{hartshorne}.
Secondly, a multi-valued map is a description
  if and only if it is homogeneous and relevant in the sense of Sections~\ref{sec_homogeneity_condition} and \ref{sect_relevance}.
One direction of this equivalence 
% that $\Phi$ is a description if and only if it is homogeneous and relevant 
is straightforward (see also \cite[Prop.~4.9]{brown_jabu_maps_of_toric_varieties}):

\begin{prop}\label{prop:description_satisfy_homog_and_relev}
        If $\Phi$ is  a description of a rational map
          $\varphi \colon X \dashrightarrow Y$, then $\Phi$ satisfies the
          homogeneity and relevance conditions.
\end{prop}

\begin{prf}
By Definition~\ref{defin:description} of description, $\pi_Y \circ \Phi$
is defined on an open subset of $\overline{X}$, so $\Phi(x)$ cannot be contained
in the irrelevant locus for those points. Therefore $\Phi$ satisfies the
relevance condition \ref{item:rel_cond_points}.

Since $\Phi$ is a description the agreement locus $\agree{\Phi}{\varphi}$
contains an open dense subset of $\Reg\Phi$.
The homogeneity condition~\ref{item:homogeneity_cond_points_on_open} is
satisfied on this set.
\end{prf}

The converse implication is slightly more involved (compare to \cite[Thm~4.10]{brown_jabu_maps_of_toric_varieties}).

\begin{thm}\label{thm:Phi_homog_and_relevant_is_a_description}
Let $\Phi\colon \overline{X} \to \overline{Y}$ be a multi-valued map that satisfies the homogeneity and relevance conditions.
Then by its action on rational functions, $\Phi^*$ naturally determines
a rational map $\varphi \colon X \dashrightarrow Y$, and $\Phi$ describes $\varphi$.
\end{thm}

\begin{prf}
   By Lemma~\ref{lem_relevant_then_R0_generates_field}
    the multivalued map $\Phi$ determines a map of fields with kernel
    $\Phi^* \colon K(Y) \dashrightarrow K(X)$.
    Let $R^0\subset K(Y)$ be the ring generating $K(Y)$ from the definition of map of fields with kernel (Definition~\ref{def_map_of_fields_with_kernel}).
    By \cite[Prop.~2.14(ii)]{brown_jabu_maps_of_toric_varieties} the ring homomorphism
    \[
      \Phi^*\colon R^0  \to K(X)
    \]
      determines a rational map $\varphi\colon X \dashrightarrow Y$ which is characterised
      by its action on rational functions $q\in K(Y)$ being $\varphi^*(q) = \Phi^*(q)$.

 Next we have to prove that $\Phi$ describes $\varphi$.
 Consider the open subset\footnote{Our choice of open subset is anticipated from Proposition~\ref{thm:agreement_locus}.}
 \[
   U=\Reg{\Phi} \setminus \Bigl(\irrel X  \ \cup \ \Phi^{-1}\bigl(\irrel Y\bigr)\Bigr).
 \]
 Note that $U$ is indeed open in $\overline{X}$ or in $\Reg \Phi$ as it is a complement of closed subsets:
   $\irrel{X}$ and $\irrel{Y}$ are closed by definition and $\Phi^{-1}(\irrel{Y})$ is closed by Proposition~\ref{prop_preimage_of_closed_is_closed}.
 Furthermore,  if $ \xi \in \Phi^{-1}(\irrel{Y})$,  or equivalently if $\Phi(\xi) \cap \irrel{Y} \ne \emptyset$,
   then $\Phi(\xi) \subset \irrel{Y}$ by \ref{item:homogeneity_cond_points}.
 So the relevance condition \ref{item:rel_cond_points} guarantees that $\Phi^{-1}(\irrel{Y}) \ne \Reg \Phi$ and U is  non-empty.
 Since $\overline{X}$ is irreducible (see \cite[Thm~1.5.1.1]{arzhantsev_derenthal_hausen_laface_Cox_rings}), the open subset is dense.
 Choose any $\xi \in U$.
 By the homogeneity condition \ref{item:homogeneity_cond_points},
 $\pi_Y(\Phi(\xi))$ is a single point $y$.
 We claim $y = \varphi(\pi_X(\xi))$,
 so that $\xi\in \agree{\Phi}{\varphi}$.

 To prove the claim, we set $x= \pi_X(\xi)$ and evaluate rational functions
 $q \in K(Y)$ at $\varphi(x)$ and~$y$:
 \[
   q(\varphi(x)) = (\varphi^* q) (x) = (\Phi^*q) (\pi_X(\xi)) = q([\Phi(\xi)]) = q(y).
 \]
 So no rational function on $Y$ can distinguish between $\varphi(x)$ and $y$ and therefore $y= \varphi(x)$.
 Hence $U \subset \agree{\Phi}{\varphi}$ and $\Phi$ describes $\varphi$.
\end{prf}

\begin{cor}\label{cor_Phi_determines_phi_and_rational_functions}
   Suppose $\varphi\colon X \dashrightarrow Y$ is a rational map of Mori Dream Spaces,
     and $\Phi \colon \overline{X} \multito \overline{Y}$ is a homogeneous and relevant multivalued map.
   Then the following conditions are equivalent:
   \begin{enumerate}
    \item \label{item_Phi_describes_phi}
          $\Phi$ describes $\varphi$.
    \item \label{item_Phi_and_phi_agree_on_a_subring}
          $\Phi^*$ and $\varphi^*$ determine the same map of fields with kernel $S(Y)^0 = K(Y) \dashrightarrow S(X)^0=K(X)$.
   \end{enumerate}
\end{cor}

\begin{prf}
   By Theorem~\ref{thm:Phi_homog_and_relevant_is_a_description} and its proof the multivalued map $\Phi$ describes
      some rational map $\psi\colon X\dashrightarrow Y$ and $\psi$ is characterised by the property that
      $\psi^* \colon K(Y) \dashrightarrow K(X)$ is equal to $\Phi^*\colon S(Y)^0 = K(Y) \dashrightarrow S(X)^0$.

   If \ref{item_Phi_describes_phi} holds, then $\psi= \varphi$
     by \cite[Lem.~I.4.1 and Def. on page 24]{hartshorne} and \ref{item_Phi_and_phi_agree_on_a_subring} is satisfied.
   Instead, assuming \ref{item_Phi_and_phi_agree_on_a_subring}, maps of fields with kernels $\varphi^*$ and $\psi^*$ are equal,
     thus also the rational maps $\varphi$ and $\psi$ are equal and \ref{item_Phi_describes_phi} holds.

\end{prf}

Now we show that every rational map has a description (see \cite[Thm~4.12]{brown_jabu_maps_of_toric_varieties}).

\begin{thm}\label{thm:existence_of_description}
   Let $\varphi\colon X \dashrightarrow Y$ be a rational map of Mori Dream Spaces.
   Then there exists a description $\Phi\colon \overline{X} \multito \overline{Y}$ of $\varphi$.
   Moreover, $\Phi$ may be chosen in such a way that it satisfies the following \emph{zeroes condition}:
   \begin{enumerate}
      \renewcommand{\theenumi}{\textnormal{(Z)}}
      \item \label{item:zeroes_cond}
            for all homogeneous regular sections $f\in S[Y]$:
            \[
              \Phi^*f = 0 \Longleftrightarrow \varphi(X) \subset \Supp (f)
            \]
            where $\Supp(f) \subset Y$ denotes the support of the divisor in $Y$ defined by $f$.
   \end{enumerate}
\end{thm}

\begin{prf}
   We must define the ring homomorphism $\Phi^*\colon S[Y] \to \overline{S(X)}$.
   We construct it in several steps.
   The construction will lead to a homomorphism from a significantly larger subring $R \subset S(Y)$,
      and in fact we are rather building a map of fields with kernel $\Phi^* \colon S(Y) \dashrightarrow \overline{S(X)}$
      in the sense of \S\ref{sec_maps_of_fields_with_kernel}.

   Note that $\Phi^*$ (considered as a map of fields with kernel) after restriction to $K(Y)\subset S(Y)$ must coincide with $\varphi^* $
      by Corollary~\ref{cor_Phi_determines_phi_and_rational_functions}.
   Moreover, let $I = I(\overline{\varphi(X)}) \subset S[Y]$ be the homogeneous ideal of (the closure of) the image of $\varphi$.
   We define the image of $f\in I$ under $\Phi^*$ to be $0$.
   Thus it remains to define a homomorphism $S[Y]/I \to \overline{S(X)}$ keeping in mind the constraints determined by $\varphi$.

   We will define the homomorphism $S[Y]/I \to \overline{S(X)}$ as a map on homogeneous elements of $S[Y]/I$
      in such a way that the multiplicative structure is preserved and that the map is linear on each degree.
   These conditions are enough to extend the map to the required ring homomorphism.
   Since $X$ is irreducible, and thus also its image $\overline{\varphi(X)}$ is irreducible.
   Thus $I$ is $G_Y$-prime, or equivalently the homogeneous elements of $S[Y]/I$ are not zero divisors.
   Therefore (as our interest are restricted to homogeneous elements),
      we are going to work with  $S[Y]/I$ as if it was a domain, in particular we will consider fractions of homogeneous elements.

   Let $G \subset \Cl(Y)$ be the weight group of $S[Y]/I$, that is the group generated by the degrees $d$ such that $(S[Y]/I)^d \ne 0$.
   Pick a finite sequence $\fromto{f_1}{f_k}$ of homogeneous sections $f_i \in S[Y]$, such that $\fromto{\deg f_1}{\deg f_k}$ generate $G$ and $f_i \notin I$.
   Define subgroups $A_i = \langle\fromto{\deg f_1}{\deg f_i}\rangle \subset G$, and define the ascending sequence of subrings:
   \[
     \kk= T_0 \subset T_1 \subset T_2 \subset \dotsb \subset T_{k} =  S[Y]/I,
   \]
   where $T_i$ comprises all gradings of $S[Y]/I$ in $A_i$.
   That is $T_i = \bigoplus_{a \in A_i} (S[Y]/I)^a$.
%    Note that if $a, b \in G$ are such that $a-b \in A_i$,
%      and $f\in S[Y]^a$, $g \in S[Y]^b$ with $g \notin I$,
%      then $\frac{f}{g}$ can be seen as

   We will gradually extend our $\kk$-algebra homomorphism $\overline{\Phi^*} \colon T_i \to \overline{S(X)}$ in such a way that for every $a\in A_i$:
   \begin{enumerate}
    \item \label{item_property_of_Ti_homogeneity} $(T_i)^a$ is mapped into the set of homogeneous multivalued sections in~$\overline{S(X)}$.
    \item \label{item_property_of_Ti_zero} if  $g\in (T_i)^a$ is mapped to $0$, then $g=0$.
    \item \label{item_property_of_Ti_agrees_with_phi} if $f,g \in (S[Y])^a$ and $g\notin I$, and $[f], [g] \in S[Y]/I$ are classes of $f$ and $g$ respectively,
             then $\frac{\overline{\Phi^*}([f])}{\overline{\Phi^*}([g])} = \varphi^*\left(\frac{f}{g}\right)$.
   \end{enumerate}
   For $i=0$, the homomorphism is defined and we proceed by induction on $i$.
   Suppose it is defined on $T_i$ for some $0 \le i <k$.
   We want to extend it it to $T_{i+1}$.
   The class $[\deg f_{i+1}]$ generates the quotient group $A_{i+1}/A_i$.
   If the group $A_{i+1}/A_i \simeq \ZZ$, then set for instance $\overline{\Phi^*}[f_{i+1}] := 1$, or any other non-zero homogeneous multivalued section in $\overline{S(X)}$.
   If $A_{i+1}/A_i$ is finite of order $r$ instead, then $[f_{i+1}]^r \in T_i$, 
      and set $\overline{\Phi^*}[f_{i+1}] := \sqrt[r]{\overline{\Phi^*}[f_{i+1}]^r}$.
   Also in this case the image of $[f_{i+1}]$ is a non-zero homogeneous multivalued section.

   Now pick any homogeneous section $g \in S[Y]^{\alpha \deg f_{i+1} + b}$ for some $b\in A_i$ and an integer $\alpha$.
   Then $g \cdot f_{i+1}^{-\alpha}\in S(Y)$ has degree $b$, hence by Lemma~\ref{lem_extending_fields_by_homogeneous_elements}
      this homogeneous rational section is in the subfield of $S(Y)$ generated by $S(Y)^0$ and $\fromto{f_1}{f_i}$.
   In other words, express $b$ as an integral combination of the generators of $A_i$, $b = \sum_{j=1}^{i} \beta_j \deg f_j$.
   Then
   \[
      q:= g \cdot f_{i+1}^{-\alpha} \cdot f_{1}^{-\beta_1} \dotsm f_{i}^{-\beta_i} \in S(Y)^0 = K(Y).
   \]
   Note the poles of $q$ are all contained in the union of zeroes of $f_i$.
   Since $f_i \notin I$ (by our choice of $f_i$), it makes sense to pullback $q$ by $\varphi$.
   Hence we define:
   \begin{equation}\label{equ_define_pullback_Phi_g}
     \overline{\Phi^*}([g]) := \varphi^*q \cdot (\overline{\Phi^*}[f_{i+1}])^{\alpha} \cdot \overline{\Phi^*} ([f_{1}]^{\beta_1} \dotsm [f_{i}]^{\beta_i})
   \end{equation}
   It is straightforward to check that the definition is linear in $g \in S[Y]^{\alpha \deg f_i + b}$ and equal to $0$
     for $g\in I$. So it does not depend on the choice of $g, g'$ in the same equivalence class $[g] = [g'] \in S[Y]/I$.
   Moreover, if $g \notin I$, then $\overline{\Phi^*}([g])$  is a homogeneous multivauled section (cf.~\ref{item_property_of_Ti_homogeneity}) and
     $\overline{\Phi^*}([g]) \ne 0$ (cf.~\ref{item_property_of_Ti_zero}).
   Also by the inductive property \ref{item_property_of_Ti_agrees_with_phi} the definition in \eqref{equ_define_pullback_Phi_g}
     does not depend on the choice of integer $\alpha$, or on the choice of expression $b = \sum_{j=1}^{i} \beta_j \deg f_j$.
   Moreover, the multiplicative structure is preserved by $\overline{\Phi^*}$ for homogeneous elements.
   Thus by linearity we extend $\overline{\Phi^*}$ to a ring homomorphism $T_{i+1} \to \overline{S(X)}$,
   which satisfies properties \ref{item_property_of_Ti_homogeneity}--\ref{item_property_of_Ti_agrees_with_phi}.

   Therefore inductively we have constructed a ring homomorphism
   \[
     \overline{\Phi^*} \colon S[Y]/I \to  \overline{S(X)},
   \]
     and by composition with the quotient, also a desired ring homomorphism $\Phi^* \colon S[Y] \to  \overline{S(X)}$.
   Note that \ref{item_property_of_Ti_homogeneity} implies that $\Phi^*$ determines a multivalued map
     $\Phi \colon \overline{X} \multito \overline{Y}$ and $\Phi$ satisfies homogeneity condition \ref{item:homogeneity_cond_all}.
   Moreover, \ref{item_property_of_Ti_zero} implies that condition \ref{item:zeroes_cond} holds,
     and hence also relevance \ref{item:rel_cond_kernel} holds.

   Thus by Theorem~\ref{thm:Phi_homog_and_relevant_is_a_description} the multivalued map is a description of a map $\psi$.
   We conclude the proof using property \ref{item_property_of_Ti_agrees_with_phi} 
     and Corollary~\ref{cor_Phi_determines_phi_and_rational_functions} which implies $\varphi = \psi$.
\end{prf}

\section{The agreement and disagreement loci}
\label{sect:agreement_revisited}

We describe the agreement locus of a description and then we prove its complement has only components of codimension $1$
  and components where the rational map is not regular.
The content of this section is parallel to \cite[\S4.4]{brown_jabu_maps_of_toric_varieties}.

\begin{prop}\label{thm:agreement_locus}
 Let $\Phi$ be a description of $\varphi$.
 Then
 \[
   \agree{\Phi}{\varphi} = \Reg{\Phi} \setminus \Bigl(\irrel X  \ \cup \ \Phi^{-1}\bigl(\irrel Y\bigr)\Bigr).
 \]
\end{prop}

\begin{prf}
 By the definition of the agreement locus,
 if $\xi \in \agree{\Phi}{\varphi}$, then
 \[
    \xi \in \Reg \Phi  \setminus \irrel X.
 \]
 The homogeneity condition holds for $\Phi$ by Proposition~\ref{prop:description_satisfy_homog_and_relev}, so, for such $\xi$,
 $\Phi(\xi)$ is contained in a single orbit
 by condition~\ref{item:homogeneity_cond_points} of
 Proposition~\ref{prop:equivalence_of_homog_conds}.
 Since $\pi_Y(\Phi(\xi))$ is defined it follows that no point in
 $\Phi(\xi)$ is in $\irrel Y$,
 which proves the first inclusion:
 \[
   \agree{\Phi}{\varphi} \subset \Reg{\Phi} \setminus
   \Bigl(\irrel X  \ \cup \ \Phi^{-1}\bigl(\irrel Y\bigr)\Bigr).
 \]

 To prove the other inclusion, take
 $\xi \in \Reg{\Phi} \setminus \Bigl(\irrel X  \ \cup
 \ \Phi^{-1}\bigl(\irrel Y\bigr)\Bigr)$
 and set $y = \pi_Y(\Phi(\xi)) \in Y$.
 We must prove, that $x=\pi_X(\xi) \in \Reg \varphi$ and that $\varphi(x) = y$,
 in other words that $\varphi^*$ maps the local ring $\ccO_{Y, y} \subset K(Y)$
 into the local ring $\ccO_{X, x} \subset K(X)$.
 So take any $q \in \ccO_{Y, y}$.
 By Corollary~\ref{cor_Phi_determines_phi_and_rational_functions},
 \[
   \varphi^*q = \Phi^*q \quad\text{as elements of $K(X)$}.
 \]
 Since a lift of $y$ to $\overline{Y}$ is in the image of $\Phi$,
 it follows that $\Phi^*q$ is defined and hence $\varphi^*q$ is defined.
 Hence we can calculate:
 \[
   (\varphi^*q)(x) = (\Phi^*q)(\xi) = q(\Phi(\xi)) = q(y),
 \]
 where the outer equalities hold because rational functions can
 be evaluated on any representative of a point in the Cox cover.
 Since $q$ is regular at $y$, also $\varphi^*q \in \ccO_{X,x}$ as claimed.
 So $\varphi(x) = y$ and thus $\xi \in \agree{\Phi}{\varphi}$.
\end{prf}

\begin{cor}\label{cor:pi_arg_is_open}
  Let $\Phi$ be a description of $\varphi$.
  The agreement locus $\agree{\Phi}{\varphi}$ is an open $G_X$-invariant
  subset of $\overline{X}$ (and of $\Reg \Phi$).
%   In addition, if $X$ is $\QQ$-factorial, then $\pi_X(\agree{\Phi}{\varphi})$ is open.
%   In general,
  $\pi_X(\agree{\Phi}{\varphi})$ contains an open dense subset of $X$.
\end{cor}

\begin{prf}
  $\Reg \Phi$ is an open $G_X$-invariant subset by
  Corollary~\ref{cor_regularity_locus_is_affine}. $\irrel X$ is clearly
  closed and $G_X$-invariant. Finally, $\irrel Y$ is a $G_Y$-invariant subset
  of $\overline{Y}$, so by homogeneity condition~\ref{item:homogeneity_cond_points}
  also $\Phi^{-1}\bigl(\irrel Y\bigr)$ is $G_X$-invariant, and it is closed in
  $\Reg \Phi$ by Proposition~\ref{prop_preimage_of_closed_is_closed}. Thus
  $\agree{\Phi}{\varphi}$ is open and $G_X$-invariant by
  Proposition~\ref{thm:agreement_locus}.
\end{prf}

The definition of the agreement locus gives
$\agree{\Phi}{\varphi} \subset {\pi_X}^{-1}(\Reg \varphi)$.
In \S\ref{sec_complete_descriptions}, we distinguish those
descriptions for which the equality holds and we prove that they always exist.
In the meantime, we  call the difference between the two sets
the \emph{disagreement locus}.

\begin{prop}\label{thm:disagreement_is_of_codimension_1}
	Let $\varphi\colon X \dashrightarrow Y$ be a rational map between
	two toric varieties $X$ and $Y$ with a description
	$\Phi\colon \overline{X} \multito \overline{Y}$.
	Consider two open subsets $U_2 \subset U_1$ of $\overline{X}$:
	\[
		U_1 = {\pi_X}^{-1}(\Reg \varphi) \quad\text{and}\quad
		U_2 = \agree{\Phi}{\varphi}.
	\]
	The disagreement locus $D=U_1 \setminus U_2$ is
	a closed subset in $U_1$ purely of codimension~$1$ in $U_1$ or is empty.
\end{prop}

\begin{prf}
   Since $U_2$ is a non-empty open subset of $U_1$
   by Corollary~\ref{cor:pi_arg_is_open},
   clearly $D$ is a proper closed subset in $U_1$.
   By Proposition~\ref{thm:agreement_locus} we have an equality
   \[
     U_2 = \Reg{\Phi} \setminus \Bigl(\irrel X  \ \cup \ \Phi^{-1}\bigl(\irrel Y\bigr)\Bigr).
   \]
   Note that $\irrel X$ is disjoint from $U_1$ (because $\pi_X$ is not regular on $\irrel X$).
   Therefore
   \begin{align*}
     D & =  \underbrace{\Bigl(U_1 \setminus  \Reg{\Phi} \Bigl)}_{=:\Dind}
            \cup
            \underbrace{\Bigl(U_1 \cap \Phi^{-1}\bigl(\irrel Y \bigr)\Bigl)}_{=:\Dirrel}.
   \end{align*}
    By Corollary~\ref{cor_regularity_locus_is_affine} the locus $\Dind$ is indeed purely of codimension~$1$ (or empty).
    It therefore remains to prove that also $\Dirrel$ is purely of codimension~$1$ or empty.

    Assume $\Dirrel$ is not empty and choose arbitrary $\xi \in \Dirrel$.
    We have to prove the codimension of $\Dirrel$ at $\xi$ is $1$.
    Since $\xi \in U_1$ the rational map $\varphi$ is regular at $x = \pi_X(\xi)$.
    Consider $y = \varphi(x)$ and its open affine neighbourhood $V=U_g \subset Y$,
    such that $V$ is given by non-vanishing of certain homogeneous regular section $g \in S[Y]$ as in Lemma~\ref{lem_open_affines_U_f}.
    Set $\gamma=\Phi^*g$.
    By homogeneity condition~\ref{item:homogeneity_cond_all}, there exists
    $f \in \kk[\Reg {\Phi}]$ such that $\gamma^r = f$ for some $r \ge 1$.
    We claim that $f(\xi)=0$ and that for all $\xi'$ in the locus $\set{f=0}$
    and in some sufficiently small open neighbourhood of $\xi$
    we have $\xi' \in \Dirrel$.

    First we prove $f(\xi) =0$.
    Since $\xi \in \Phi^{-1}\bigl(\irrel Y \bigr)$
    it follows $\Phi(\xi)$ and $\irrel Y$ have non-empty intersection.
    As usual, since $\Phi(\xi)$ is contained in a single $G_Y$-orbit
    by the homogeneity condition~\ref{item:homogeneity_cond_points},
    we have $\Phi(\xi) \subset \irrel Y$.
    In particular, $\Phi(\xi)$ is disjoint from ${\pi_Y}^{-1} (V)$,
    in other words the section $g$ vanishes on $\Phi(\xi)$.
    So $\gamma$ vanishes at $\xi$ and therefore $f$ vanishes on $\xi$.

    We prove further that $\xi' \in \set{f=0}$ implies $\xi' \in \Dirrel$,
    at least in some neighbourhood of $\xi$.
    More precisely, we take this neighbourhood to be
    \[
       (\varphi \circ {\pi_X})^{-1}(V) \cap \Reg \Phi.
    \]
    Since $\Phi$ is regular at such $\xi'$:
    \[
      0= f(\xi') = \gamma^r(\xi') = \Phi^*(g^r) (\xi') = g^r (\Phi(\xi')),
    \]
    so $\Phi(\xi')$ is contained in the locus $g =0$.
    Therefore $\Phi(\xi')$ is disjoint from ${\pi_Y}^{-1}(V)$
    and hence the set $\pi_Y(\Phi(\xi'))$ (if non-empty) is not in $V$.
    On the other hand $\varphi(\pi_X(\xi'))$ is contained in $V$ by our choice of the open
    neighbourhood of $\xi$.
    We conclude, that $\xi'$ cannot be in the agreement locus $U_2$.
    But $\xi' \in \Reg{\Phi}$ and $\xi' \notin \irrel X$ (again by our choice
    of open neighbourhood of $\xi$).
    Therefore by Proposition~\ref{thm:agreement_locus} there is no other
    possibility than $\xi' \in \Phi^{-1}\bigl(\irrel Y\bigr)$ so that
    $\xi' \in \Dirrel$ as claimed.

    Hence $\Dirrel$  locally near $\xi$ contains a subset $\set{f=0}$ purely of
    codimension~$1$.
    Since the same holds true for every $\xi\in \Dirrel$ and $\Dirrel \ne U_1$,
    we conclude that $\Dirrel$ is purely of codimension~$1$.
\end{prf}

\section{Existence of complete descriptions}
\label{sec_complete_descriptions}

\renewcommand{\theenumi}{\textnormal{(\Alph{enumi})}}
\renewcommand{\labelenumi}{\theenumi}

\begin{defin}
\label{defin:complete_agreement}
A description $\Phi$ of $\varphi\colon X \dashrightarrow Y$ is \emph{complete}
if it satisfies
 \begin{enumerate}
   \setcounter{enumi}{2}
   \item \label{item:complete_agreement}
           $\agree{\Phi}{\varphi} = \pi_X^{-1} (\Reg \varphi)$.
 \end{enumerate}
\end{defin}

Proposition~\ref{thm:agreement_locus}, together with this definition,
has an immediate corollary.

\begin{cor}\label{cor:regularity_locus_from_complete_description}
  If $\Phi$ is a complete description of $\varphi$, then
  \[
    \Reg \varphi = \pi_X( \Reg \Phi  \setminus \Phi^{-1}(\irrel Y )).
  \]
  In particular, the map $\varphi$ is regular on $X$ if and only if
  $\Phi$ is regular on $\overline{X}$ and $\Phi^{-1}(\irrel Y)$ is contained in $\irrel X$.
\end{cor}

The main claim of this article is that complete descriptions of maps between Mori Dream Spaces always exist.
\begin{thm}\label{thm:existence_of_complete_descriptions}
  Let $\varphi\colon X \dashrightarrow Y$ be a rational map of Mori Dream Spaces.
  Then there exists a complete description $\Phi\colon \overline{X} \multito \overline{Y}$ of $\varphi$.
   Moreover, $\Phi$ may be chosen in such a way that it satisfies zeroes condition \ref{item:zeroes_cond}
      of Theorem~\ref{thm:existence_of_description}.
\end{thm}
The theorem is proved throughout this section.
We commence with an additional assumption that $Y$ is a toric variety.
For this part, we follow the strategy of the arguments of \cite[\S4.5]{brown_jabu_maps_of_toric_varieties},
   where both $X$ and $Y$ are assumed to be toric varieties.
Later, to conclude we use the embedding theorem for Mori Dream Spaces.
Every Mori Dream Space admits a natural closed embedding into a toric variety, see \S\ref{sec_general_case_of_existence_of_complete_descriptions} for details and references.

\subsection{Case when the target is a toric variety}
  Suppose $Y$ is a normal toric variety, with Cox ring $S[Y] \simeq \kk[\fromto{y_1}{y_n}]$,
    where the generators $y_i$ are homogeneous of degrees $\deg(y_i) \in \Cl(Y)$,
    and they correspond to rays of the fan $\Sigma_Y$ (see \cite[\S5.1]{cox_book}).
    
  Assume $w \colon \Cl(Y) \to \QQ$ is a group homomorphism. 
  If the image is contained in $\ZZ$, then we can define the action of $\kk^*$ on $\overline{Y}$ as follows.
  Suppose $\eta \in \overline{Y} \simeq \kk^n$ is a closed point $\eta = (\fromto{\eta_1}{\eta_n})$ and pick $t \in \kk^*$.
  Then define $t^w\cdot \eta := (\fromto{t^{w(\deg{y_1})} \eta_1}{t^{w(\deg{y_n})} \eta_n})$.
  To extend the definition to rational characters $w$, 
    we must be careful to choose the correct roots, similarly to~\cite[\S3.4 and \S1.1.2]{brown_jabu_maps_of_toric_varieties}.
  If the image of $w$ is not contained in $\ZZ$, then there exists an integer $r$ such that $r\cdot w \colon \Cl(Y) \to \ZZ $.
  Let $s = \sqrt[r]{t}$ be any root, and then define:
  \[
      t^w\cdot \eta :=  (\fromto{s^{r\cdot w(\deg{y_1})} \eta_1}{s^{r\cdot w(\deg{y_n})} \eta_n})
  \]
  That is, we pick a root $s$ once and for all, and we apply the same root on all coordinates of $\overline{Y}$. 
  With this definition, the following lemma is clear.
  \begin{lemma}\label{lem_action_on_points_of_Y}
     Suppose $\eta \in \overline{Y} \simeq \kk^n$ is a point $\eta = (\fromto{\eta_1}{\eta_n})$.
     Pick a group homomorphism $w \colon \Cl(Y) \to \QQ$ and any $t \in \kk^*$.
     Then the point $t^w\cdot \eta$  is in the same $G_Y$-orbit as $\eta$
       (independent of the choice of the root $s$ as above).
  \end{lemma}
  \noprf
  Now we add to the picture a Mori Dream Space $X$ and a rational map $\varphi\colon X \dashrightarrow Y$.
  By Theorem~\ref{thm:existence_of_description} there exists a description
  \[
    \Phi\colon \overline{X} \multito \overline{Y} \simeq \kk^n
  \]
  satisfying zeroes condition \ref{item:zeroes_cond}.
  What are the possibilities to modify $\Phi$, so that it still describes $\varphi$ and satisfies \ref{item:zeroes_cond}?

  If $f \in S[X]$ is a homogeneous regular section and $w\colon \Cl(Y) \to \QQ$ is a group homomorphism,
     then we define a multi-valued map
\begin{align*}
    f^w \cdot \Phi\colon \overline{X} &\multito \overline{Y}\\
    \xi  &  \multimapsto f(\xi)^w\cdot \Phi(\xi) = \left(\fromto{f^{w(\deg{y_1})}\Phi^*y_1}{f^{w(\deg{y_n})}\Phi^*y_n}\right).
\end{align*}
\begin{cor}\label{cor_agreement_locus_of_modification}
    $f^w\cdot\Phi$ describes $\varphi$ and their agreement loci are the same away from $\set{f=0}$:
    \[
       \agree{\Phi}{\varphi} \cap \set{f\ne 0} =  \agree{f^{w} \cdot \Phi}{\varphi} \cap \set{f\ne 0}.
    \]
    Moreover, $f^w\cdot\Phi$ satisfies \ref{item:zeroes_cond}.
\end{cor}
\begin{prf}
  On a sufficiently general point $\xi$ of $\overline{X}$ (such that $f(\xi) \ne 0$) the image $(f^w \cdot \Phi)(\xi)$
    is in the same $G_Y$-orbit as the image of $\Phi(\xi)$ by Lemma~\ref{lem_action_on_points_of_Y}.
  Thus  $f^w\cdot\Phi$ describes $\varphi$ and the agreement locus is as claimed.
  Furthermore, say $g \in S[Y]$ is a homogeneous section of degree $\deg g$.
  Then
  \[
     (f^w \cdot \Phi)^* g =  f^{w(\deg g)} \Phi^*g
  \]
  and since $f^{w(\deg g)}$ is invertible in $\overline{S(X)}$ we have
  \[
     (f^w \cdot \Phi)^* g = 0 \Longleftrightarrow \Phi^* g = 0 \Longleftrightarrow \varphi(X) \subset \Supp(g).
  \]
  Thus the condition \ref{item:zeroes_cond} is satisfied.
\end{prf}

We want to consider the toric variety $Z \subset Y$, which is the smallest (closed) toric stratum of $Y$,
  which contains $\varphi(X)$.

Let $\sigma \in \Sigma_Y\subset N\otimes\RR$ be the cone corresponding to $Z$,
and let $\Sigma_{Y,Z}=\{\tau_Y \in\Sigma_Y \mid \sigma \preceq \tau_Y\}$
%%% THINK: formally, we must include also faces of such $\tau_Y$, so that $\Sigma_{Y,Z}$ is a fan. But probably, this is sufficiently clear.
be the subfan of $\Sigma_Y$
corresponding to the smallest torus invariant open neighbourhood in $Y$
 of the toric stratum $Z$. Let $N_\sigma=\sigma\cap N$ be the sublattice of $N$ generated by $\sigma$ and let $N(\sigma)=N/N_\sigma$
 be the quotient lattice. Let $\Sigma_Z=\starfan{\sigma}$ be the fan given by the images of cones in $\Sigma_{Y,Z}$ under the natural projection
    $p \colon N\longrightarrow N(\sigma)$; this is the fan of $Z$ regarded as a toric variety under the $T_{N(\sigma)}$ torus action, cf.~\cite[\S 3.2]{cox_book}
(if $\varphi(X)$ is not contained in any strict toric stratum
of $Y$, then $Z=Y$ and $\quofan$ is equal to $\Sigma_Y$).

Denote the primitive  generators of rays of $\Sigma_Y$ by $\fromto{\rho_1}{\rho_n}$,
   where the variable $y_i$ in the Cox ring $S[Y]$ corresponds to the ray spanned by $\rho_i$. Let $\widehat{N}=\bigoplus_{i=1}^{n}\ZZ e_{\rho_i}$
be the group of torus invariant Weil divisors on $Y$. Define $L$ to be the natural surjective linear map $L:\widehat{N}\longrightarrow N(\sigma)$ sending $e_{\rho_i}$ to $p(\rho_i)$ and set $L_\QQ=L\otimes \QQ$.
The kernel of $L_\QQ$ is responsible for the freedom we have in modifying the descriptions.
\begin{lemma}\label{lem_modify_description}
   Suppose $f\in S[X]$ is a non-zero homogeneous section and  $\mu=(\mu_1,\ldots,\mu_n)\in\ker L_\QQ$.
   Then the multi-valued map $\Psi \colon \overline{X} \multito \overline{Y}$ given by
   \[
     \Psi^*y_i := f^{\mu_i} \Phi^*y_i.
   \]
   is also a description of $\varphi$.
   Moreover $\Psi$ satisfies zeroes condition~\ref{item:zeroes_cond} and
     \[
       \agree{\Phi}{\varphi} \cap \set{f\ne 0} =  \agree{\Psi}{\varphi} \cap \set{f\ne 0}.
    \]
\end{lemma}
\begin{prf}
    Set $R=\{i\mid \rho_i\in\sigma(1)\}$ and $R'=\{1,\ldots,n\}\setminus R$. The zeroes condition~\ref{item:zeroes_cond} and the choice of $Z\subset Y$ assure that $\Phi^* y_i=0$ if and only if $i\in R$.
    The vector space $N_{\sigma}\otimes \QQ$ is spanned by $\set{\rho_i : i \in R}$ and, by the assumption, $\sum_{j\in R'}\mu_j \rho_j\in N_\sigma \otimes \QQ$.
    Choose $\nu_i\in\QQ$ such that the latter is equal to $\sum_{i\in R}\nu_i\rho_i.$
    Set $\mu'=\sum_{j\in R'}\mu_j e_{\rho_j}-\sum_{i\in R}\nu_i e_{\rho_i}$ and  $\mu''=\sum_{i\in R}(\nu_i+\mu_i) e_{\rho_i}$.
    Then $\mu=\mu'+\mu''$ and $\mu'$ defines a homomorphism $w\colon\Cl(Y) \to \QQ$ by setting $w(\deg y_i) = -\nu_i$
       for $i\in R$ and $w(\deg y_j)=\mu_j$ for $j\in R'.$

   Note $f^{\mu_i} \Phi^* y_i = f^{-\nu_i}\Phi^*y_i =0$ for $i \in R$.
   Therefore
   \[
      f^{\mu_i} \Phi^*y_i = f^{w(e_i)} \Phi^* y_i,
   \]
   or equivalently, $\Psi = f^w\cdot \Phi$ and the claim follows
     from Corollary~\ref{cor_agreement_locus_of_modification}
\end{prf}

Fix a prime divisor $(f) \subset X$ (in particular, $f\in S[X]$ is homogeneous and $G_X$-irreducible).
For all $i \in \setfromto{1}{n}$ let $\mu_i$ be the multiplicity of $f$ in $\Phi^* y_i$ if the pullback is non-zero, or pick any number $\mu_i\in \QQ$, if $\Phi^* y_i=0$. Set $\mu=\sum_{i=1}^n \mu_i e_{\rho_i}\in \widehat{N}\otimes\QQ.$

\begin{lemma}\label{lem:vanishing_order_along_divisor}
Let $m\in\Hom_\ZZ(N(\sigma),\ZZ)$ be an integral linear form corresponding to the rational function $\chi^m$ on $Z$.
Then, for $f\in S[X]$ as above, $L_\QQ(\mu)$ is an integral point of $N(\sigma)$ and the order of vanishing of  $\varphi^* \chi^m$
     along the divisor $(f)$ is equal to $m\circ L(\mu)$.
\end{lemma}

  \begin{prf}
  The composition of $m\circ p$ is an integral form on the lattice $N$ and hence $\chi^{m\circ p}$ is a  rational function on $Y$ whose
  order of vanishing along the toric divisor corresponding to the ray spanned by $\rho_i$ is equal to $m\circ L(e_{\rho_i})$.
  The order of $\Phi^* y_i$ along $(f)$ is by definition $\mu_i$
    so  the order of $\varphi^*\chi^m$  (which is equal to $\Phi^* \chi^{m\circ p}$ by Corollary~\ref{cor_Phi_determines_phi_and_rational_functions})
  along $(f)$ is $m \circ L_\QQ(\mu)$. Since $m$ is arbitrary $L_\QQ(\mu)\in N(\sigma)$ is integral.
  \end{prf}

\begin{cor}\label{cor:Lv_in_support_of_fan}
    For $f\in S[X]$ as above, if $L_\QQ(\mu)$ is not in the support of $\quofan$, then $\varphi$ is
       not regular on $(f)$.
\end{cor}

\begin{prf}
    Let $\tau$ be any cone in $\quofan$.
    Since $L(\mu) \notin \tau$, there exists $m_{\tau}\in \Hom_\ZZ(N(\sigma),\ZZ)$ in the dual cone $\tau^{\vee}$ such that
       $m_{\tau}\circ L (\mu) < 0$.
    Then, by Lemma~\ref{lem:vanishing_order_along_divisor},
    the rational function $\varphi^* \chi^{m_{\tau}}$ has a pole along $(f)$.
    Let $U_{\tau}$ be the affine open subset of $Y$ corresponding to a cone in $\Sigma_{Y,Z}$,
    which maps onto $\tau$.
    Note that the collection of such $U_{\tau}$ for all $\tau \in \quofan$
       will cover the image of $\varphi$.
    Then $m_\tau\circ p$ is a regular function on $U_\tau$ for which the pullback is not regular on $(f)$. By \cite[Prop.~2.15]{brown_jabu_maps_of_toric_varieties}
       this implies that $\varphi$ is not regular on $(f)$.
\end{prf}

We are now ready to prove:
\begin{prop}\label{prop:existence_of_complete_descriptions_target_toric}
  Let $X$ be a Mori Dream Space and $Y$ a toric variety.
  Suppose $\varphi\colon X \dashrightarrow Y$ is a rational map.
  Then there exists a complete description $\Phi\colon \overline{X} \multito \overline{Y}$ of $\varphi$.
  Moreover, $\Phi$ may be chosen in such a way it satisfies zeroes condition \ref{item:zeroes_cond}
      of Theorem~\ref{thm:existence_of_description}.
\end{prop}

\begin{prf}
By Theorem~\ref{thm:existence_of_description} there exists a description $\Phi$ of $\varphi$ satisfying \ref{item:zeroes_cond}.
      We have to modify this description to obtain a complete one.

By Proposition~\ref{thm:disagreement_is_of_codimension_1},
the disagreement locus
\[
    D = \pi_X^{-1} (\Reg \varphi) \setminus \agree{\Phi}{\varphi}
\]
is a union of codimension $1$ components.
If $D$ is empty, then the proposition is proved, so suppose it is not empty; we must modify $\Phi$
so that the new description is defined on those components which cover the locus where $\varphi$ is defined.

Choose any homogeneously prime component of $D$ and pick a
   $G_X$-ir\-re\-du\-cible section $f\in S[X]$ that vanishes to order $1$ along this component.
As above, for all $i \in \setfromto{1}{n}$
   let $\mu_i$ be the multiplicity of $f$ in $\Phi^* y_i$ if the pullback is non-zero,
  or pick any number $\mu_i\in \QQ$, if $\Phi^* y_i=0$.

By Corollary~\ref{cor:Lv_in_support_of_fan}, if $L(\mu)$ does not lie in the support of $\quofan$, then $(f)$
is not a part of the disagreement locus, contradicting our setup.
Thus  $L(\mu)$ lies in the support of $\quofan$.

Let $\tau$ be the cone in $\quofan$ of minimal dimension that contains $L(\mu)$
  and let $\tau_Y$ be a maximal cone in $\Sigma_{Y,Z}$ that maps onto $\tau$.
Pick a vector $u\in \tau_Y$ that maps to $L(\mu)$.
There exists a vector $\mu'=(\mu'_1,\ldots,\mu'_n)\in\widehat{N}$ with
$\mu'_i>0$ for $\rho_i\in\tau_Y(1)$ and $\mu'_i=0$ otherwise, such that $p(\mu')=u$.
Since $\mu'-\mu \in \ker L$, by Lemma~\ref{lem_modify_description}, the multi-valued map $\Psi$
\[
    \Psi^*y_i:= f^{\mu'_i - \mu_i} \cdot \Phi^*y_i
\]
is another description of $\varphi$ with the same (dis)agreement locus as $\Phi$ away from $\set{f=0}$.
We claim that the component $\set{f=0}$ is not in the disagreement
  of $\Psi$ and $\varphi$.

By Proposition~\ref{thm:agreement_locus}, it is enough to prove the
following two statements:
\begin{itemize}
    \item $\Psi$ is regular on a general point of $(f)$.
    \item $\Psi$ does not map a general point of $(f)$
	into the irrelevant locus of $Y$.
\end{itemize}
The first is immediate: $\Psi$ is regular along $(f)$ since the multiplicity of $f$ in any non-zero $\Psi^* y_i$ is equal to $\mu_i'$ which is non-negative.
Moreover, this shows that if $x \in \set{f=0}$ is a general point, then $\Psi(x)$ has zero $y_i$-coordinate
if and only if either $\Phi^*y_i=0 $ or $\mu'_i > 0$.
In particular, if $\rho_i\notin\tau_Y(1)$,
  then the $i$-th coordinate of $\Psi(x)$ is non-zero.
This means that the generator of $B_Y$ determined by $\tau_Y$ is non-zero
  at $\Psi(x)$ so $x$ is not mapped by $\Psi$ into the irrelevant locus of $Y$.
Therefore $\agree{\Psi}{\varphi}$
contains a general point of $\set{f=0}$ as claimed.

Thus we have obtained a description $\Psi$ of
$\varphi$ whose disagreement locus contains one component less
than that of $\Phi$. Continuing inductively, we obtain a description
with an empty disagreement locus, namely a complete description.
\end{prf}

\subsection{General case}\label{sec_general_case_of_existence_of_complete_descriptions}

Now let $Y$ be a Mori Dream Space.
Then there exists a toric variety $Z$ and a closed embedding $\iota\colon Y \to Z$,
  such that $\Cl(Z) \simeq \Cl(Y)$  and the Cox rings satisfy the graded equality $S[Y] = S[Z]/I$,
  where $I = I(\iota(Z)) \subset S[Z]$ is the homogeneous ideal of sections vanishing on $\iota(Y)$,
  see \cite[Thm~3.2.1.4, Construction~3.2.5.3, Prop.~3.2.5.4(i),(iii))]{arzhantsev_derenthal_hausen_laface_Cox_rings}.
Using this embedding we can reduce the general case of Theorem~\ref{thm:existence_of_complete_descriptions} to
the case when the target is a toric variety.

\begin{prf}[ of Theorem~\ref{thm:existence_of_complete_descriptions}]
  Suppose $X$ and $Y$ are Mori Dream Spaces, and $\varphi\colon X\dashrightarrow Y$ is a rational map.
  Let $Z$ be the toric variety and  $\iota\colon Y \to Z$ be the closed embedding as above.
  Set $\psi:=\iota\circ \varphi$ to be the composed rational map.
  By Proposition~\ref{prop:existence_of_complete_descriptions_target_toric},
    there exists a complete description $\Psi\colon \overline{X} \multito \overline{Z}$ of $\psi$,
    satisfying in addition zeroes condition \ref{item:zeroes_cond}.
  In particular, for all sections $f \in I = I(\iota(Y))$, the pullback $\Psi^*f$ vanishes, because $\psi(X) \subset \iota(Y)$.
  Thus $\Psi^*$ factorises through $S[Z]/I = S[Y]$:
  \[
     S[Z] \to S[Z]/I = S[Y] \stackrel{\Phi^*}{\to} \overline{S(X)}.
  \]
  We define the multi-valued map $\Phi \colon \overline{X} \multito \overline{Y}$
    by the algebra homomorphism $\Phi^*$ arising from this factorisation.
  To make sure this is well defined, we must check that on a distinguished set of homogeneous generators $y_i\in S[Y]$
    the image $\Phi^*y_i$ is homogeneous. But homogeneity condition~\ref{item:homogeneity_cond_all} for $\Psi$
    guarantees in fact a stronger statement, that the same homogeneity condition holds for $\Phi$.
  Note that here we exploit that $\Cl(Y) = \Cl(Z)$, and thus $G_Y = G_Z$.

  We claim $\Phi$ is a complete description of $\varphi$.
  Firstly note that $\Reg \Phi = \Reg \Psi$,
    since the denominators of images of $\Phi^*\colon S[Y] \to \overline{S(X)}$ are the same as the denominators of the images of $\Psi^*\colon S[Z] \to \overline{S(X)}$
    (simply, the two images are equal).
  Pick $\xi \in \widehat{X}$ such that $\varphi$ is defined at $\pi_X(\xi)$.
  Then also $\psi = \iota \circ \varphi$ is defined at $\pi_X(\xi)$ and thus (since $\Psi$ is a complete description) both $\Psi$ and $\Phi$ are regular at $\xi$.
  Moreover, $\Psi(\xi)$ is contained in a single $G_Z$-orbit, the sections of $I$ (treated as functions on $\overline{Z}$) vanish on $\Psi(\xi)$,
     and the image $\pi_Z (\Psi(\xi))$ is a single point in $Z$ equal to $\psi\circ\pi_X(\xi) = \iota \circ \varphi\circ\pi_X(\xi)$.
  Therefore the $G_Z$ orbit containing $\Psi(\xi)$ is contained in $\overline{Y} \subset \overline{Z}$, and $\Phi(\xi)$ is in the same $G_Y$-orbit (again, recall $G_Y= G_Z$).
  Thus $\pi_Y(\Phi(\xi))$ is a single point and it is equal to $ \varphi\circ\pi_X(\xi)$ as claimed.
\end{prf}

%  \bibliography{references}
\def\polhk#1{\setbox0=\hbox{#1}%
  {\ooalign{\hidewidth\lower1.5ex\hbox{`}\hidewidth\crcr\unhbox0}}}\def\dbar{\leavevmode\hbox
  to 0pt{\hskip.2ex\accent"16\hss}d}

\end{document}